\newif\ifenglish\englishtrue
\ifenglish
\catcode127=11
\fi
\def\probel#1{\ifx#1.\relax\else\ \fi#1}
\hsize18truecm
\hoffset-1truecm
\vsize25truecm
\voffset-.5truecm
\frenchspacing
%
\font\small=cmr8             %
\font\ss=cmss10              %
\font\ss=cmss10              %
\font\ssdbf=cmssdc10         %
\font\ssqi=cmssqi8           %
\font\vtt=cmvtt10            %
          %
              %
             %
          %
             %
%
\def\hexnumber#1{\ifcase#1 0\or 1\or 2\or 3\or 4\or 5\or 6\or 7\or 8\or
                 9\or A\or B\or C\or D\or E\or F\fi}
\newfam\relfam
\font\tenmsam=msam10
\font\sevenmsam=msam7
\font\fivemsam=msam5
\textfont\relfam\tenmsam
\scriptfont\relfam\sevenmsam
\scriptscriptfont\relfam\fivemsam
\mathchardef\le"3\hexnumber\relfam36
\mathchardef\ge"3\hexnumber\relfam3E
\mathchardef\preceq"3\hexnumber\relfam34
\mathchardef\succeq"3\hexnumber\relfam3C


\def\R{\hbox{\MMMM R}}

\font\ZV=cmr10.pk scaled 2000
\def\zvezda{\mathop{\,\vrule width0pt depth2pt height8pt
            \smash{\lower7pt\hbox{\ZV *}}\,}\limits}
\font\angles=cmex10 scaled 2000
\newfam\anglesfam
\textfont\anglesfam=\angles

\def\mapsfrom{\mathrel{\leftarrow\kern-1pt\vrule height4.2pt depth-.3pt}}
\newfam\trafam
\font\tenmsbm=msbm10
\font\sevenmsbm=msbm7
\font\fivemsbm=msbm5
\textfont\trafam\tenmsbm
\scriptfont\trafam\sevenmsbm
\scriptscriptfont\trafam\fivemsbm
\mathchardef\varnothing="0\hexnumber\trafam3F
\let\emptyset\varnothing
%
\let\epsilon\varepsilon
\let\phi\varphi
\let\~\widetilde
\let\=\overline
\let\^\widehat

\let \< \langle
\let \> \rangle

\def\:#1{\ifmmode
           \ifx#1={\buildrel {\rm\ifenglish def\else ®¯à\fi} \over =}
           \else\colon #1
           \fi
         \else{\rm :}#1%
         \fi
        }
\def\-{\item{-}}
\def\s#1\par{\message{#1}\par\goodbreak\bigskip\leftline{\bf #1}%
        \nobreak\par\smallskip\noindent%
}
\def\gp#1{\left\langle#1\right\rangle}
\def\nc#1{\gp{\!\gp{#1}\!}}

\def\1{\{1\}}
\def\0{\{0\}}
%

\ifenglish\else
\csname newcount\endcsname\«­®
\def\‹\par{
\par
\advance\«­® by 1
\goodbreak\vtop to 0pt{\llap{\it ‹.\the\«­®\ }\hfill\vss\nobreak}
\nointerlineskip
\par
}
\csname newcount\endcsname\nozad
\def\§{%
  \advance\nozad by 1%
  \¡ãª­®=159%
  \ifhmode\par%
        \else\smallskip%
  \fi%
  \goodbreak%
 \item{\bf\the\nozad.}
}
\csname newcount\endcsname\¡ãª­®
\def\¡ãª{
  \ifnum\¡ãª­®<166\char\¡ãª­®
  \else\ifnum\¡ãª­®=166 ñ\else \advance \¡ãª­® by -1 \char\¡ãª­® \advance \¡ãª­® by 1
              \fi
  \fi
}
\¡ãª­®=159
\def\¯{
  \advance\¡ãª­® by 1
  \itemitem{\¡ãª)}
}
\def\"ª § ­¨¥{{\it "ª § ­¨¥.\ }}
\let\"\"ª § ­¨¥
\def\latc#1{%
\ifnum`#1>127
\if#1 a\fi
\if#1¡b\fi
\if#1¢c\fi
\if#1£d\fi
\if#1¤e\fi
\if#1¥f\fi
\if#1¦g\fi
\if#1§h\fi
\if#1¨i\fi
\if#1ªj\fi
\if#1«k\fi
\if#1¬l\fi
\if#1­m\fi
\if#1®n\fi
\if#1¯o\fi
\if#1àp\fi
\if#1áq\fi
\if#1âr\fi
\if#1ãs\fi
\if#1ät\fi
\if#1åu\fi
\if#1æv\fi
\if#1çw\fi
\if#1èx\fi
\if#1éy\fi
\if#1êz\fi
\else#1\fi
}
\def\latr#1{\ifx\relax#1\else\latc#1\expandafter\latr\fi}
\def\"®ª § â¥«ìá⢮{\noindent{\bf "®ª § â¥«ìá⢮.\ }}
\let\"\"®ª § â¥«ìá⢮
\def\Ž¯à¥¤¥«¥­¨¥{\noindent{\bf Ž¯à¥¤¥«¥­¨¥.\ }}
\let\Ž¯à\Ž¯à¥¤¥«¥­¨¥
\def\'¥®à¥¬ {\csname proclaim\endcsname '¥®à¥¬ \probel}
\let\'\'¥®à¥¬ 
\def\‹¥¬¬ {\csname proclaim\endcsname ‹¥¬¬ \probel}
\def\"⢥ত¥­¨¥{\csname proclaim\endcsname "⢥ত¥­¨¥\probel}
\let\"â¢\"⢥ত¥­¨¥
\def\à¥¤«®¦¥­¨¥{\csname proclaim\endcsname à¥¤«®¦¥­¨¥\probel}
\def\à¨¬¥à#1. #2\par{\csname proclaim\endcsname
        à¨¬¥à\probel#1. {\rm #2}\par}
\def\‡ ¬¥ç ­¨¥#1. #2\par{\csname proclaim\endcsname
    ‡ ¬¥ç ­¨¥\probel #1. {\rm #2}\par}
\def\ƒ¨¯®â¥§ {\csname proclaim\endcsname ƒ¨¯®â¥§ \probel}
\def\'«¥¤á⢨¥{\csname proclaim\endcsname '«¥¤á⢨¥\probel}
\fi
\ifenglish
\def\Proof{\noindent{\bf Proof.\ }}

\def\Theorem{\csname proclaim\endcsname Theorem\probel}
\let\Th\Theorem
\def\Lemma{\csname proclaim\endcsname Lemma\probel}
\def\Proposition{\csname proclaim\endcsname Proposition\probel}
\def\Example#1. #2\par{\csname proclaim\endcsname
        Example\probel#1. {\rm #2}\par}

\def\Remark#1. #2\par{\csname proclaim\endcsname
        Remark\probel#1. {\rm #2}\par}

\def\Conjecture{\csname proclaim\endcsname Conjecture\probel}
\def\Corollary{\csname proclaim\endcsname Corollary\probel}
\fi
\def\({{\rm(}}
\def\){{\rm)}}
\ifenglish\else
\hyphenation{£àã¯-¯ë}
\fi
\catcode127=15
\hfuzz5pt
\def\A{{\cal A}}
\def\R{{\cal R}}
\def\T{{\cal T}}
\def\X{{\cal X}}
\def\Y{{\cal Y}}

\count0=1
\headline{\ifnum\count0=1 \vtt UDC 512.543.7+512.546.1\hss\else\hss\fi}

\centerline{\ssdbf THE NUMBER OF NON-SOLUTIONS TO AN EQUATION IN A GROUP}
\centerline{\ssdbf
AND NON-TOPOLOGIZABLE TORSION-FREE GROUPS%
\footnote{*}{\rm This work was supported by the Russian Foundation for
Basic Research, project no. \number02-01-00170.}}

\smallskip
\centerline{\ss Anton A. Klyachko and 
Anton V. Trofimov}
\smallskip
{
\ssqi
\centerline{Faculty of Mechanics and Mathematics, Moscow State University}
\centerline{Moscow 119992, Leninskie gory, MSU}
\centerline{klyachko@daniil.math.msu.su}
\centerline{anton\_tr@rambler.ru}
}
\medskip
{\narrower\small\noindent
It is shown that, for any pair of cardinals with infinite sum, there exist
a group and an equation over this group such that the first cardinal is the
number of solutions to this equation and the second cardinal is the number of
non-solutions to this equation. A countable torsion-free non-topologizable
group is constructed.

\noindent{\it Key words\/\rm:}
equations in groups, topologization of groups, graded presentations.

\noindent{\it AMS MSC\/\rm:} 20F05, 20F06, 22A05, 54H11.

}
\medskip

\s 1. Introduction

\proclaim{Main theorem}.
There exists a finitely generated torsion-free group $H$ and an
equation $w(x)=1$ over $H$ \(this means that $w(x)\in H*\gp x_\infty$\)
such that the set of solutions to this equation consists of
all elements of $H$ exept one:
$$
\{h\in H\ |\ w(h)\ne1\}=\1.
\eqno{(1)}
$$

\Remark.
The L\"owenheim--Skolem theorem implies that the words ``finitely
generated" in the main theorem can be replaced by the words ``of arbitraty
infinite cardinality" (if (1) is valid in a group $H$, then it is
also valid in, e.g., all ultrapowers of $H$ and in all their
subgroups containing the diagonal).

Since the set of solutions to an equation is closed in any Hausdorff
group topology, we obtain the following corollary, which answers 
a question of P.~I.~Kirku [UPTA85, Question 1.4].

\Corollary.
There exists a nontrivial countable torsion-free group admitting no
non-discrete Hausdorff group topology.

Note that, according to Markov's theorem [M46], the complement to an
element in any countable non-topologizable group decomposes into a
finite union of the sets of solutions to some systems of equations.
In the known examples of infinite countable non-topologizable groups,
these decompositions look as follows:
$$
\eqalign{
G\setminus\{g_1,\dots,g_n\}=\bigcup_{i=1}^{n-1}\{g\in G\ |\ g^n=g_i\}
\quad&\hbox{(Ol$'$shanskii's example$^{**}$ 
[O80], [O89] and its
modifications [MO98])};\cr
G\setminus\{g_1,\dots,g_{2n}\}=\{g\in G\ |\ [g,a]^n=1\}
\quad&\hbox{(examples from [T04])}.\cr
}
$$
\footnote{}{\llap{$^{**}$ }Ol$'$shanskii's example is a central quotient 
of a group constructed by Adyan [A71].} 
Here $g_i$ and $a$ are some fixed elements of the corresponding
group $G$, and the number $n$ in both cases is large (at least 665)
and odd. Decomposition (1) is simplest possible.
However, the equation $w(x)=1$ constructed in the proof of the
main theorem is much more complicated (see formulae $(**)$
and $(*)$) than the equations $x^n=a$ and $[x,a]^n=1$ from the
earlier known examples of countable non-topologizable groups. Note also
that Ol$'$shanskii's example and its modifications are periodic
(the examples of Morris and Obraztsov [MO98] are also quasicyclic); the
examples from [T04] have torsion, but are not periodic; moreover, in
[T04], it was shown that any countable group can be embedded in one of
such examples.

It is natural to ask, what values the cardinality of the
set of solutions to an equation in a group and the cardinality
of the complement to this setcan take. The following fact can be
easily derived from the main theorem.

\Th{1}.
For any two cardinals $\bf s$ and $\bf n$ at least one of which is
infinite, there exists a group $G$ \(of cardinality ${\bf s}+\bf n$\) and
an equation $u(x)=1$ over $G$ such that precisely $\bf s$ elements of $G$
are solutions to this equation and precisely $\bf n$ elements of $G$ are
not solutions to this equation.

\Remark.
The condition about infinity of one of the cardinals cannot be omitted.
For example, it is easy to see that, in the group of order three, no
equation has exactly two solutions.

The authors thank A. Yu. Ol$'$shanskii for reading this paper and useful 
comments and corrections. 


\s 2. Ol$'$shanskii's approach to constructing groups with prescribed
properties

A {\it graded presentation} is a group presentation
$G(\infty)=\gp{\A\ |\ \R}$ with some filtration on the set of relators:
$$
\R=\bigcup_{i=0}^\infty\R_i,\quad
\emptyset=\R_0\subseteq\R_1\subseteq\dots.
$$
Relators from $\R_i\setminus\R_{i-1}$ are called {\it relators of rank
$i$}, and the presentation $\gp{\A\ |\ \R_i}$ is denoted by $G(i)$ and
considered as a graded presentation (assuming $\R_j=\R_i$ for $j>i$).

According to [O89], we say that a graded presentation
$G(\infty)$ without periodic relations satisfies the {\it condition} R
(with parameters $\alpha$, $h$, $d$, and $n$) if, for each $i$, there
exists a set $\X_i\subseteq F=\gp{\A\ |\ \varnothing}$ of words (called
the {\it periods of rank~$i$}) such that

\item{1)}
          the length of each word from $\X_i$ is $i$ and
          no word from $\X_i$ is conjugate in $G(i-1)$
          to a power of a word of smaller length;

\item{2)}
          if $X$ and $Y$ are different words from $\X_i$, then $X$ is not
          conjugate to $Y^{\pm1}$ in $G(i-1)$;

\item{3)}
          each relator $R\in\R_i\setminus\R_{i-1}$ has the form
$$
R\equiv\prod_{k=1}^h(T_kA^{n_k}),\hbox{ where } A\in\X_i,
$$
and the following conditions hold:

\itemitem{R1.} $n_k\ge n$;

\itemitem{R2.} $|n_i|/|n_j|\le 1+{1\over2}h^{-1}$;

\itemitem{R3.} The words $T_k$ are not equal in $G(i-1)$ to words of
               smaller length, and $|T_k|<di$;

\itemitem{R4.} $T_k\notin\gp A$ in the group $G(i-1)$;

\itemitem{R5.} The word $R$ is not a proper power in the free group,
               and if
               $V\equiv A^{n_{s-1}}\prod\limits_{k=s}^{s+l}(T_kA^{n_k})$
               is a cyclic subword of $R$, $l\ge \alpha^{-1}-4$, and
               $VV_1$ and $VV_2$ are cyclic shifts of $R$, then
               $V_1\equiv V_2$;

\itemitem{R6.} Let $V\equiv A^{m_1}T_kA^{n_k}\dots T_{k+l}A^{m_2}$ be
               a subword of a cyclic shift of $R$, where
               $l\ge\alpha^{-1}-2$, and let
               $V'\equiv A^{m_1'}T_{k'}'A^{n_k'}\dots T_{k'+l}'A^{m_2'}$
               be a subword of a cyclic shift of a relator
               $(R')^{\pm1}$ with the same period $A$; suppose that the
               signs of the exponents $m_1$ and $m_1'$, $n_k$ and $n_k'$,
               \dots, $m_2$ and $m_2'$ coincide. If $V$ and
               $V'$ decompose into concatenations of words: $V\equiv
               V_0\dots V_l$ and $V'\equiv V_0'\dots V_l'$, where
               $V_0\equiv A^{m_1}T_kA^{c_1}$, $V_1\equiv
               A^{b_1}T_{k+1}A^{c_2}$, \dots, $V_l\equiv
               A^{b_l}T_{k+l}A^{m_2}$, $V_0'\equiv
               A^{m_1'}T_{k'}'A^{c_1'}$, $V_1'\equiv
               A^{b_1'}T_{k'+1}'A^{c_2'}$, \dots, $V_l\equiv
               A^{b_l'}T_{k'+l}'A^{m_2'}$, and $V_j=V_j'$ for $j=0,\dots,l$
               in the group $G(i-1)$,
               then
               $R'\equiv R$, $V'\equiv V$, and $V$ is not a subword
               of the cyclic word $R^{-1}$.

\noindent
In the book [O89], one can find many useful properties of presentations
with condition~R. We need some of these properties.

\Lemma 1.
If a graded presentation $G(\infty)$ satisfy condition
{\rm R} with sufficiently small number $\alpha$ and sufficiently
large numbers $h$, $d$, and $n$ \($1\ll\alpha^{-1}\ll h\ll d\ll n$\), then
\item{\rm1)}
         all abelian subgroup of the group $G(\infty)$ are cyclic;
\item{\rm2)}
         the group $G(\infty)$ is torsion-free;
\item{\rm3)}
         if elements $X$ and $Y$ are conjugate in $G(\infty)$, then there
         exists an element $Z\in G(\infty)$ such that $X=ZYZ^{-1}$ and
         $|Z|\le({1\over2}+\alpha)(|X|+|Y|)$;
\item{\rm4)}
         if $A$, $B$, and $C$ are nonidentity elements of the group
         $G(\infty)$ and $X$ is an element such that $X^{-1}AXB$ and $C$
         are conjugate in $G(\infty)$, then the double coset $\gp A X\gp B$
         contains an element $X'$ of length smaller than
         $({1\over2}+\alpha)(|A|+|B|+|C|)+
         \left[{1\over2}|A|\right]+\left[{1\over2}|B|\right]$
         \rm(here the square brackets denote the integer part of a number);
\item{\rm5)}
         if a word $X$ represents the identity element of $G(\infty)$,
         then $X=1$ in the group $\gp{\A\ \bigm|\ \{R\in\R\ |\
         |R|<(1-\alpha)^{-1}|X|\}}$;
\item{\rm6)}
         if elements $X$ and
         $ZXZ^{-1}$ commute in $G(\infty)$, then $X$ and $Z$ commute in
         $G(\infty)$.

\Proof
All these properties were proved in [O89]:
the first property is one of the assertions of Theorem 26.5;
the second property is a special case of Lemma 25.2;
the third property is Lemma 25.4;
the fourth property is a translation into the algebraic language of
a slightly weakened assertion of Lemma 22.2 about cuts of diagrams on
a sphere with three holes;
the fifth property is a translation into the algebraic
language of Lemma 23.16;
and the sixth property coincides with Lemma 25.14.


\s 3. Construction of the group $H$

Take a sufficiently large even number $h$ and an integer
$n\gg h$. As the alphabet $\A$, we take the set of letters
$\{a,b,c_1,c_2,\dots,c_h\}$. Let $F$ be the free group with basis $\A$.

We set
$\R_0=\R_1=\R_2=\emptyset$. Let us define the presentation $G(i)=\gp{\A\ |\
\R_i}$ for $i>2$ assuming that the presentation $G(i-1)=\gp{\A\ |\
\R_{i-1}}$ is already defined. As the set of periods of rank $i$,
we take some set $\X_i\subset F$ of words of length~$i$ satisfying
conditions

\item{1)}
          no word from $\X_i$ is conjugate in $G(i-1)$
          to a word of smaller length;

\item{2)}
          different words from $\X_i$ are not conjugate in $G(i-1)$ to
          each other and to the inverses of each other;

\item{3)}
          each word from $\X_i$ equals $ab$ in
          $G(i-1)/([G(i-1),G(i-1)]\gp{c_1c_2\dots c_h})$;

\item{4)}
          the set $\X_i$ is maximal among all sets
          satisfying conditions 1), 2), and 3).

\noindent
Let us define a set $\Y_{i,j,Z}\subseteq F$, where $j=1,\dots,h$ and
$Z\in G(i-1)$,
as the set of all minimal (i.e., not equal to words of smaller length)
in $G(i-1)$ words of length smaller than $di$ representing
the element $Zc_jZ^{-1}$ in $G(i-1)$. Finally, we set
$$
\R_i=\R_{i-1}\cup\T_i,\quad
$$
where $\T_i$ is a maximal set of pairwise nonconjugate in $G(i-1)$ words
of the form
$$
R=R_{A,Z}=\prod_{j=1}^h\left(T_jA^{(-1)^jn}\right),\quad\hbox{ where }
A\in \X_i,\ T_j\in\Y_{i,j,Z},\ Z\in G(i-1).
$$

\Lemma 2.
The graded presentation
$$
H=G(\infty)=\gp{\A\ |\ \R}, \quad\hbox{where }\ \R=\bigcup_{i=0}^\infty\R_i,
$$
satisfies the condition {\rm R}.

\Proof
Note that, modulo the commutator subgroup, we have:
\item{-}
     each relation of the presentation $H$ has the form $c_1c_2\dots c_h=1$;

\item{-}
     $T_j=c_j$;

\item{-}
     each period has the form $ab(c_1\dots c_h)^k$.

\noindent
Therefore, the periods are not proper powers and
conditions 1), 2), R1, R2, and R3 hold by the construction.
Condition R4 obviously holds, even modulo the commutator
subgroup.

It is well-known that a long common subword of two
words of the form $A^n$ must be consistent; hence, the violation of
condition R5 implies that either one of the words $T_j$ lies in $\gp
A$ or two words $T_j$ with different subscripts $j$ lie in the same
double coset with respect to $\gp A$. But neither of these two
possibilities is realized, even modulo the commutator
subgroup.

Suppose that condition R6 does not hold. We have an equality of the form
$T_j^{\pm1}A^p=A^sT_{j'}'$ in $G(i-1)$. Considering this equality modulo
the commutator subgroup, we conclude that $j=j'$, $p=s$, and
$\pm1=1$. Recalling that $T_j=Zc_jZ^{-1}$ and $T_j'=Z'c_j(Z')^{-1}$ in
$G(i-1)$, we see that $c_j$ commutes with $Z^{-1}A^pZ'$ in $G(i-1)$.
By the induction, we can assume that the presentation $G(i-1)$
satisfies condition R, in particular, commuting elements of $G(i-1)$
must lie in the same cyclic subgroup (by Lemma 1); since $c_j$ is not
a proper power (even modulo the commutator subgroup), we obtain
an equality of the form $Z^{-1}A^pZ'=c_j^k$. But the violation of
condition R6 for $G(i)$ means that similar equalities hold in
$G(i-1)$ for many numbers $j$. Considering two such equalities
$$
Z^{-1}A^pZ'=c_j^k\quad\hbox{and}\quad Z^{-1}A^{p_1}Z'=c_{j_1}^{k_1}
$$
(with $j\ne j_1$) modulo the commutator subgroup, we obtain
that $k=0$, i.e., $Z'\in\gp A Z$ and the words $R$ and $R'$ represent
conjugate elements of the group $G(i-1)$. By construction, this means that
$R\equiv R'$.


\s 4. Proof of the main theorem

\Lemma 3.
For each element $g$ of the commutator subgroup of $H$ and each word $A$
conjugate in $H$ to the element $g^{-1}agb$, there exists a
word $Z$ of length at most ${1\over3}d|A|$ such that
$Z^{-1}AZ=g^{-1}agb$ in $H$.

\Proof
According to assertion 4) of Lemma 1, $|a^pgb^r|<|A|+2$ for some integers
$p$ and $r$ (here and in what follows, $|\cdot|$ denotes the length of an
elements in $H=G(\infty)$). But considering the element $a^pgb^r$ modulo
the commutator subgroup, we see that $|a^pgb^r|\ge|p|+|r|$ and, therefore,
$$
|g|\le|p|+|r|+|a^pgb^r|\le2|a^pgb^r|<2|A|+4\quad\hbox{and}\quad
|g^{-1}agb|<4|A|+10.
$$
The last inequality and assertion 3) of
Lemma 1 implies that there exists an element $Z$ such that
$Z^{-1}AZ=g^{-1}agb$ and
$$
|Z|\le|A|+|g^{-1}agb|<5|A|+10<{1\over3}d|A|
$$
(the last inequality holds because $d\gg1$).

Consider the equation
$$
v(x)=1,\quad\hbox{where }
v(x)\equiv \prod_{j=1}^h\left(c_j(x^{-1}axb)^{(-1)^jn}\right)
\in H*\gp x_\infty,
\eqno{(*)}
$$
over the group $H$.

\Lemma 4.
In the group $H$, any nonidentity element of
the commutator subgroup is a solution to equation $(*)$, while the
identity element is not a solution to this equation.

\Proof
The inequality $v(1)\ne1$ follows from assertion 5) of Lemma 1, because
$v(1)\ne1$ in the free group,
the length of each relator $R$ of the group $H$ is at least
$3(n-d-2)h$, and
$$
{|v(1)|\over|R|}\le{2n+1\over3(n-d-2)}<1-\alpha
\quad\hbox{for } \alpha<{1\over3}\hbox{ and } n\gg d.
$$

Now, let us show that any nonidentity element $g$ of the commutator
subgroup of $H$ is a solution to equation~$(*)$. Let $u\in F$ be a word
of minimal length representing an element conjugate to $g^{-1}agb$.

If $|u|\le2$, then either $u=ab$ or $u=ba$ (it suffices to consider $u$
modulo the commutator subgroup). According to assertion 4) of Lemma 1, $u$
is conjugate to an element of the form $\~g^{-1}a\~gb$, where
$\~g\in\gp a g\gp b$ and $|\~g|<3$. Considering $\~g$ modulo the commutator
subgroup, we see that
$$
\~g
\hbox{ is}
\hbox{ either }1,\
a^{\pm1},\
b^{\pm1},\
a^{\pm1}b^{\pm1},\
\hbox{ or }b^{\pm1}a^{\pm1}.
$$
The first 4 cases are impossible, because they imply that
$g\in \gp a\gp b\cap [H,H]=\1$. To verify that the equality
$\~g=b^{\pm1}a^{\pm1}$ is impossible, it is sufficient to consider
the relation $\~g^{-1}a\~gb\sim u\sim ab$ modulo the normal subgroup
$\nc{c_1,\dots,c_h}$ generated by $\{c_1,\dots,c_h\}$, since the group
$G(\infty)/\nc{c_1,\dots,c_h}$ is free with basis $\{a,b\}$.

We have shown that $|u|>2$. But in this case, $u$ is conjugate to
a period $A$ of rank $|u|$ (by the definition of the set $\X_{|u|}$).
By Lemma 3, $Z^{-1}AZ=g^{-1}agb$ for some word $Z$ of length at most
${1\over3}d|A|$. Therefore, for each $j=1,\dots,h$, some minimal word
$T_j$ representing the element $Zc_jZ^{-1}$ in $G(|u|-1)$ has length
smaller than $d|A|$, and hence, by the construction of $H$, the word
$$
\prod_{j=1}^h\left(T_jA^{(-1)^jn}\right)
$$
is conjugate to one of the relators of $G(|u|)$. This completes the proof
of the lemma.

\Lemma 5.
The set of solutions to the equation
$$
[c_1v([a,x])c_1^{-1},v([b,x])]=1
\eqno{(**)}
$$
over $H$ is $H\setminus\1$.

\Proof
Since $a$ is not a proper power (modulo the commutator subgroup),
Lemma 1 implies $[a,g]\ne1$ for $g\notin\gp a$. Therefore,
by Lemma 4, all elements of $H$ not lying in $\gp a$
are solutions to $(**)$. For the same reasons, all elements of $H$ not
lying in $\gp b$ are solutions to $(**)$. But $\gp a \cap\gp
b =\1$, which can be easily verified by considering again the quotient
group by the commutator subgroup.

It remains to prove that the identity element is not a solution to $(**)$.
Assuming the converse, we have \break$[c_1v(1)c_1^{-1},v(1)]=1$. By
assertion 6) of Lemma 1, this implies that $[c_1,v(1)]=1$. This equality,
in its turn, means that $v(1)\in\gp{c_1}$ (by Lemma 1 and because $c_1$ is
not a proper power (modulo the commutator subgroup)). Considering the
equality $v(1)=c_1^k$ modulo the commutator subgroup, we come to
the conclusion that $k=0$ and $v(1)=1$, which contradicts Lemma 4.

Lemma 5 is proven, and the main theorem is also proven,
because, according to Lemma 1, the group $H$ is torsion-free.


\s 5. Proof of Theorem 1

If either ${\bf n}=0$, ${\bf s}=0$, or ${\bf s}=1$, then the ``equations" 
$1=1$, $g=1$ (where $g$ is a nontrivial element of the group), and $x=1$ 
over a group of suitable cardinality obviously have the required 
properties.

If $1<{\bf s}\le \bf n$, then as $G$ we can take the free
product of any abelian group $A$ of cardinality $\bf s$ and
any group $B$ of cardinality $\bf n$. It is easy to see that
the set of solutions to the equation $xa=ax$, where $a\in A\setminus\1$,
is $A$. Thus, the number of solutions is $\bf s$, and the number
of non-solutions is $\bf n$ (because $\bf n$ is infinite
and not smaller than $\bf s$).

If ${\bf s}>{\bf n}>0$, then as a $G$ we can take the direct
product of any group $K$ of cardinality $\bf n$ and a group $H$ of
cardinality $\bf s$, whose existence is asserted by the main theorem and
the remark after it. The equation $w(x)=1$ from the main theorem has
the required property.
Indeed,
equation $(**)$ constructed in the proof of the main theorem has zero
exponent sum of $x$; hence, the set of solutions to this equation in the
group $G=H\times K$ is $(H\setminus\1)\times K$, and the set of
non-solutions is $\1\times K$. This implies the assertion of the theorem.


\s{\rm REFERENCES}

\item{[A71]}
S. I. Adyan,
{\it On some torsion-free groups},
Izvest. Acad. Nauk SSSR, ser. Mat.,
{\bf 35} (1971), no. 3, pp. 459--468.

\item{[M46]}
A. A. Markov,
{\it On unconditionally closed sets},
Mat. Sbornik
{\bf 18} (1946), no. 1., pp. 3--28.

\item{[MO98]}
S. A. Morris and V. N. Obraztsov,
{\it Nondiscrete topological groups with many discrete subgroups},
Topology Appl.
{\bf 84} (1998), pp. 105--120.

\item{[UPTA85]}
{\it Unsolved problems of topological algebra},
(eds. V.~I.~Arnautov, A.~V.~Arkhangel$'$skii, P.~I.~Kirku, 
A.~V.~Mikhalev, Yu.~N.~Mukhin, I.~V.~Protasov, M.~M.~Choban),
Shtiintsa, Kishinev, 1985.

\item{[O80]}
A. Yu. Ol$'$shanskii,
{\it A remark on a countable nontopologized group},
Vestnik Moskov. Univ., Ser. I Mat., Mekh.
no.3. (1980), p. 103.

\item{[O89]}
A. Yu. Ol$'$shanskii,
{\it Geometry of defining relations in groups},
Nauka, Moscow, 1989;\hfil\break
English translation in Math. And Its Applications (Soviet series),
{\bf 70}, Kluwer Acad. Publishers, Dordrecht, 1991.

\item{[T04]}
A. V. Trofimov,
{\it An embedding theorem for non-topologizable group},
Vestn. Moskov. Univ., Ser. I,  Mat., Mekh.
(to appear).

\end